%% file: ex_article.tex
\documentclass[hidelinks,onefignum,onetabnum]{siamart251216}

\input{ex_shared}

\ifpdf
\hypersetup{
  pdftitle={On necessary conditions for feedback stabilizability when the control space is a fiber bundle},
  pdfauthor={W. Jongeneel, M. Scolamiero, and K. H. Johansson}
}
\fi

\newcommand{\mylabel}[2]{%
  \hypertarget{#1}{#2}%
}

\newcommand{\myref}[2]{%
  \hyperlink{#1}{#2}%
}

\usepackage{tikz-cd}
\usetikzlibrary{arrows}

\newcommand\xqed[1]{%
  \leavevmode\unskip\penalty9999 \hbox{}\nobreak\hfill
  \quad\hbox{#1}}
\newcommand\exampleEnd{\xqed{$\circ$}}

\usepackage{enumerate}
\usepackage{amssymb}

\begin{document}

\maketitle

\begin{abstract}
Within control theory, a question of great interest is to understand when a stabilizing feedback exists. Building upon Krasnosel'ski\u{\i} and co-workers, necessary conditions for the existence of a continuous stabilizing feedback have been pioneered by Brockett and Coron, providing so-called adversarial and homological conditions, respectively. Kvalheim recently generalized their conditions from equilibria to compact attractors and---in line with earlier work by Coron---showed that when a control system can be identified with a vector bundle, the homology condition is stronger, that is, the homology condition implies the adversary condition. Kvalheim also showed that this comparison breaks when extending the case of vector bundles to fiber bundles. In this work, we propose and characterize an appropriate adversarial condition for control systems represented by fiber bundles. The homology condition, on the other hand, is already well-defined, showcasing its general applicability. With the generalized conditions at hand, we show that by appropriately relaxing the regularity of the feedback under consideration, that is, by allowing for mild discontinuities, the homology condition implies the adversarial condition for fiber bundles. Thereby, in this work, we extend these earlier comparative results.  
\end{abstract}



\section{Introduction}
Within the broad field of automatic control, it is of great importance to understand, \textit{a priori}, if a stabilization problem admits a solution. For instance, suppose we are given a control system $\dot{x}=u$, with $x\in \mathbb{R}$ the \textit{state} and $u\in \mathbb{R}$ the \textit{input}. As $u$ is unconstrained, for all states, the so-called \textit{control space} (see \S\ref{sec:control:systems}) can be identified with the trivial vector bundle $\pi_1:X\times U\to X$, defined through $(x,u)\mapsto\pi_1(x,u):=x$, where $X:=\mathbb{R}$, $U:=\mathbb{R}$ and for any $x\in X$ the fiber is simply $\pi_1^{-1}(x)=U$, see \S\ref{sec:nec:cond:1.2} for more on bundles. The fiber captures the space of available inputs, given a state $x$. Now, suppose that the goal is to \textit{stabilize} (see \S\ref{sec:Lyapunov stability}) the origin using continuous feedback, that is, finding a continuous map $\mu_1:\mathbb{R}\to \mathbb{R}$ (more precisely, a section $X\to X\times U$, see \S\ref{sec:control:systems}) such that the origin is stable under $\dot{x}=\mu_1(x)$. Then, as the control system is linear and controllable, such a feedback is guaranteed to exist. One may simply select $x\mapsto \mu_1(x):=-x$. What if, in addition, the constraint $|u|=1$ must hold? Satisfying this constraint using linear feedback is impossible; the control space fails to be a vector bundle and is effectively the fiber bundle $\pi_2:X\times \mathbb{S}^0\to X$, with non-contractible fibers $\pi_2^{-1}(x)=\mathbb{S}^0\cong \{-1,1\}$. Here, we employ standard notation and recall that ``$\cong$'' denotes a context-dependent isomorphism. To continue, the feedback $x\mapsto \mu_2(x):=-\mathrm{sgn}(x)$ \textit{does} lead to stabilization of the origin. Compared to utilizing $\mu_1$, the key difference is that the closed-loop vector field is not continuous anymore. On the other hand, the \textit{solutions} of the corresponding differential equation, in the sense of Carath\'eodory, Krasovskii and Filippov, are continuous, forward in time. Can all of this---especially when working with more intricate control systems---be understood in a principled manner? That is the central theme of this note.

What we are looking for is a better understanding of structural features of a control system that enable or preclude stabilization by feedback of a certain regularity. With the focus on uncovering structural features, we study necessary conditions for solving feedback stabilization problems. 

To be more explicit, we illustrate the conceptual framework by highlighting the most influential results in this area. First, Brockett's seminal result \cite[Thm. 1.(iii)]{ref:brockett1983asymptotic}. Consider a control system $\dot{x}=f(x,u)$, with $f:\mathbb{R}^n\times \mathbb{R}^m\to \mathbb{R}^n$ being continuously differentiable\footnote{Brockett used tools from differential topology to prove his result and it is well-known that $C^1$ regularity is not needed for statements like the Poincar\'e-Hopf theorem to remain true, \textit{e.g.}, typically $C^0$ suffices instead of $C^1$ or even $C^{\infty}$ smoothness \cite{ref:milnor65}. Indeed, using the index theoretic results from \cite[\S 52]{ref:krasnosel1984geometrical}, Zabczyk generalized Brockett's condition to the $C^0$ setting \cite[\S 2]{ref:zabczyk1989}.} ($C^1$) and $f(0,0)=0$. Then, \textbf{\textit{Brockett's condition}} states that there is a $C^1$ feedback $x\mapsto \mu(x)$ such that $0$ is locally asymptotically stable under the closed-loop vector field $x\mapsto f(x,\mu(x))$, only if $(x,u)\mapsto f(x,u)$ maps a neighbourhood of $(0,0)$ onto a neighbourhood of $0$. Topological necessary conditions of this form have been highly influential as one can draw relatively strong conclusions based on generally available information, \textit{e.g.}, see Example~\ref{ex:nonholonomic} for a system that is controllable, but cannot be stabilized using continuous feedback, as follows from applying Brockett's condition. This should be contrasted with linear controllable systems where controllability implies the existence of a globally stabilizing continuous (linear) feedback \cite[Ch. 2]{ref:wonham1974linear}.      

Another influential necessary condition is due to Coron \cite[Thm. 2]{ref:Coron1990}. Here, we consider the same control systems as Brockett. Let $\Sigma := \mathrm{dom}(f)\setminus f^{-1}(0)$, that is, $\Sigma$ is the subset of the control space that cannot be mapped to $0$ under $f$. Then, \textit{\textbf{Coron's condition}} states that\footnote{For $n=1$, \eqref{equ:Coron:condition} extends, but using reduced homology. We already point out that Mansouri's generalization, to compact codimension 1 attractors, also relies on $n>1$, \textit{e.g.}, \cite[Thm. 4]{ref:mansouri2007local} fails to hold for $\dot{x}=-x$ with $x\in \mathbb{R}$.}, for $n\geq 2$, there is a $C^0$ feedback $x\mapsto \mu(x)$ such that $0$ is locally asymptotically stable under the closed-loop vector field $x\mapsto f(x,\mu(x))$, only if
\begin{align}
\label{equ:Coron:condition}
    f_* H_{n-1}(\Sigma) = H_{n-1}(\mathbb{R}^n\setminus \{0\})\cong \mathbb{Z}.
\end{align}
Here, $H_k(\cdot)$ denotes the $k$-th singular homology group with integer coefficients and $f_*$ denotes the corresponding induced homomorphism, \textit{e.g.}, consult \cite[Ch. IV]{ref:bredon1993topology} or \cite[Ch. 2]{Hatcher} for more on the relevant algebraic topology. Coron's condition can be shown to be stronger than Brockett's condition \cite[p. 229]{ref:Coron1990}, and as with Brockett, this condition is implied by the index theoretic results stated in \cite[\S 52]{ref:krasnosel1984geometrical}, see also~\cite{ref:bobylev1974deformation,ref:krasnosel1968operator}. On the other hand, Brockett's condition is arguably easier to work with. 

Besides relaxing differentiability to continuity, compared to original work by Brockett, generalizations of Brockett's condition pertain to stabilization of rather generic compact sets (not just equilibrium points) where the control system is defined using structures reminiscent of fiber bundles (not just an Euclidean product space) \cite{ref:kvalheim2021necessary}. The same is true for generalizations of Coron's condition \cite{ref:kvalheim2022obstructions,ref:kvalheim2023relationships}. 

Then, recently, Kvalheim showed, in great generality, how well-known necessary conditions for continuous feedback stabilization compare \cite{ref:kvalheim2023relationships}. To be precise, one may focus on 
(i) the \textit{adversary condition}, that is, Brockett's condition \cite{ref:brockett1983asymptotic} and its generalizations \cite{ref:ryan1994brockett,ref:orsi2003necessary,ref:byrnes2008brockett,ref:kvalheim2021necessary}; and (ii) the \textit{homology condition}, that is, Coron's condition \cite{ref:Coron1990} and its generalizations \cite{ref:mansouri2007local,ref:mansouri2010topological,ref:kvalheim2022obstructions,ref:Mansouri2025}. See also \cite{ref:wilson1967structure,ref:sontag1980remarks,ref:krasnosel1984geometrical,ref:zabczyk1989,ref:ceragioli2002some,ref:coron2007control,ref:moulay2011conley,ref:whydowetwist2025,ref:kvalheim2026stabilizability} for related work in this space. Here, it is interesting to note that most of these results essentially hinge on the homotopy invariance of the topological degree of a map (\textit{e.g.}, see \cite[p. 81, 108]{ref:guillemin2010differential}), that is, one typically computes some invariant for a canonical system and then relates this invariant to a homotopy class of the canonical system, \textit{e.g.}, see \cite[Ex. 3.4]{ref:jongeneel2023topological}. Now, regarding the comparison itself, Kvalheim shows that when a control space (\textit{i.e.}, the domain of $f$) can be understood as a vector bundle with oriented base space, then, the homology condition is stronger than the adversary condition, that is, if the homology condition holds, so does the adversary condition \cite[Thm. 4]{ref:kvalheim2023relationships}; the result of Coron extends. Surprisingly, the situation for fiber bundles is less clear, \textit{e.g.}, when working with sphere bundles, there are counterexamples by Kvalheim~\cite[Ex. 4]{ref:kvalheim2023relationships} and de Sa and Ames~\cite[Ex. 4]{ref:de2026topological} that show cases where the homology condition holds but the adversarial condition does not (see also Example~\ref{ex:counter} below). One might ask the following. 
\\\\
\textbf{Question A:} ``\textit{What do these counterexamples convey if one is solely interested in stabilization? Is the homology condition truly occasionally weaker than the adversarial condition?}" 
\\\\
An affirmative answer to the latter would be surprising given prior results by Coron and Kvalheim. The goal of this work is to resolve this conundrum. More details regarding Question A are provided below, in particular regarding the meaning of the word ``\textit{truly}'' as those counterexamples are evidently correct within the context of those papers. 

\textbf{Contribution}
We show in this note that the aforementioned counterexamples are \textit{not} truly counterexamples if one is solely interested in stabilization. The reason that for these examples the homology condition holds, while the adversarial condition does not, is simply the result of the homology condition being necessary for a larger class of feedback controllers; feedback that is possibly discontinuous. Indeed, those counterexamples admit a stabilizing feedback, just not a continuous one. 

\begin{figure}[th!]
    \centering
    \includegraphics[width=0.95\linewidth]{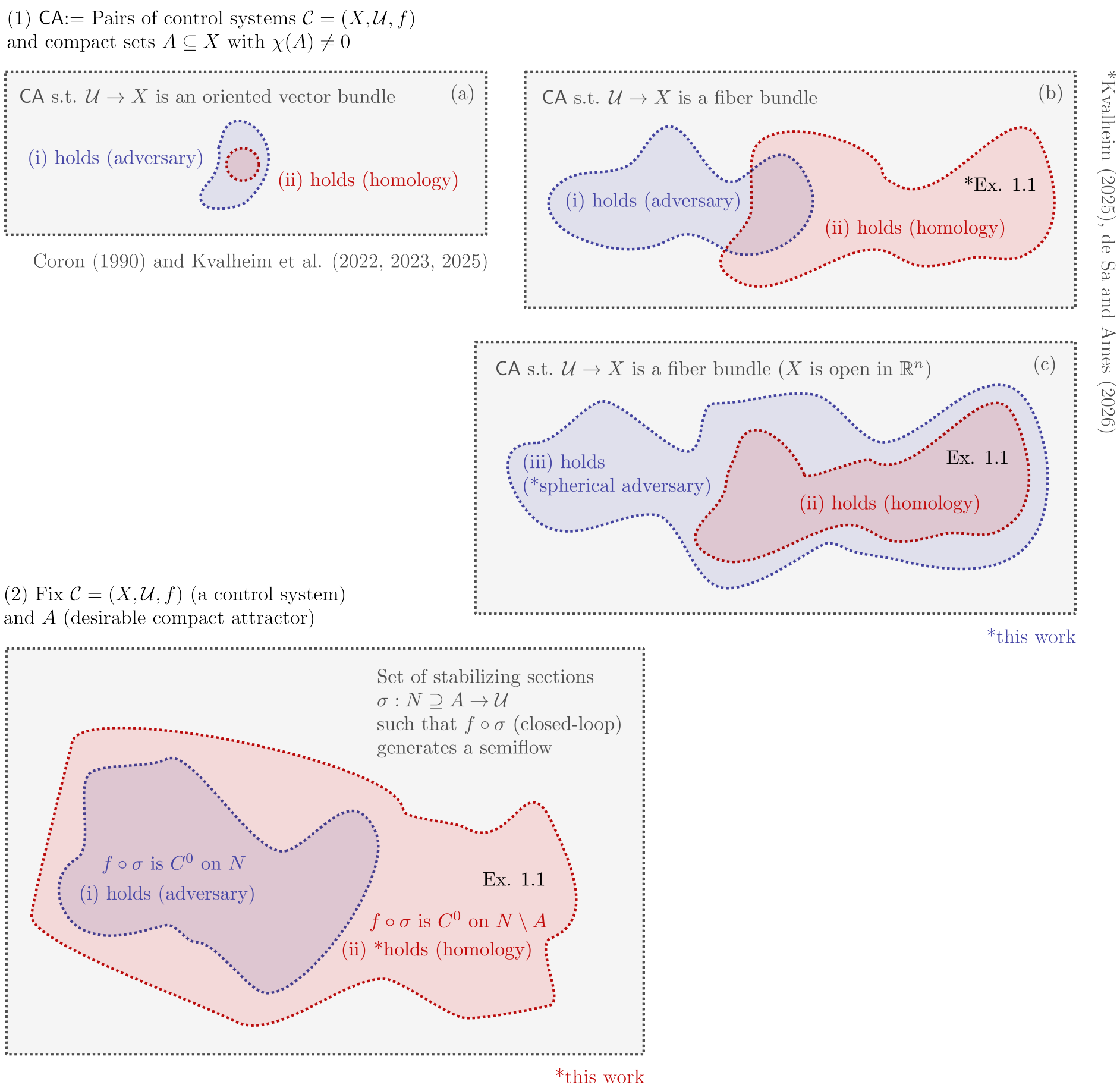}
    \caption{A summary of the context and contribution of this work. $(1a)$ When the control system (see \S~\ref{sec:control:systems}) is such that $\mathcal{U}\to X$ is a vector bundle, with $X$ oriented, then, the homology condition {\protect\myref{homology:condition}{(ii)}} is stronger than the adversary condition~{\protect\myref{adversary:condition}{(i)}}. $(1b)$ When $\mathcal{U}\to X$ is a fiber bundle (\textit{e.g.}, a sphere bundle), there are examples showcasing that this ordering fails (see Example~\ref{ex:counter}). Nevertheless, the system studied in Example~\ref{ex:counter} does admit a stabilizing controller, however, one that is merely continuous on $N\setminus A$ (not on all of $N$). $(2)$ Indeed, we show that the homology condition {\protect\myref{homology:condition}{(ii)}} holds true for stabilizing control systems with mild discontinuities on $A$, hence one could argue that $(1b)$ displays shortcomings of the adversary condition, not of the homology condition. $(1c)$ To overcome these shortcomings, we develop the spherical adversary condition {\protect\myref{spherical:adversary:condition}{(iii)}}; this adversarial condition accommodates mild discontinuities, handles Example~\ref{ex:counter} and is again implied by the homology condition.}
    \label{fig:comparison}
\end{figure}

We first show that the case of fiber bundles \textit{can} be better understood by relaxing the regularity of the feedback. In particular: (i) we construct a new adversarial condition (see Theorem~\ref{thm:spherical:adversary}); (ii) show that this condition clarifies the fiber bundle counterexample by Kvalheim (see Examples~\ref{ex:counter}-\ref{ex:finite:time:sphere}); and (iii) we show that the homology condition remains stronger (see Corollary~\ref{cor:homology} and Theorem~\ref{thm:compare}). Hence, the contribution of this note is to provide clarity beyond vector bundles and concurrently we provide further motivation for using the homology condition. We have summarized the above in Figure~\ref{fig:comparison}. The goal of this work is not to be as general as possible (\textit{e.g.}, we do not work with generic fiber bundles). Lastly, we point to open problems and conclude in Section~\ref{sec:discussion}.

\subsection{Control systems} 
\label{sec:control:systems}
To study further necessary conditions for feedback stabilization, we start by defining a (continuous) \textit{\textbf{control system}} via the triple $\mathcal{C}:=(X,\mathcal{U},f)$. Here, $X$ is the state space (an $n$-dimensional smooth ($C^{\infty}$) manifold), $\mathcal{U}$ a topological space together with a continuous surjective map $\pi_u:\mathcal{U}\to X$ (\textit{e.g.}, a fiber bundle, see \S\ref{sec:nec:cond:1.2}) known as the \textbf{\textit{control space}} and lastly, $f:\mathcal{U}\to TX$ is a continuous fiber-preserving map defined on the control space, that is, given a continuous section $\sigma:X\to \mathcal{U}$ (\textit{i.e.}, $\pi_u\circ \sigma = \mathrm{id}_X$), $f\circ \sigma:X\to TX$ is a vector field. As $\sigma$ can be understood as \textbf{\textit{feedback}}, $f\circ \sigma$ should be understood as a \textit{closed-loop} vector field, \textit{e.g.}, consider the following commutative diagram:
\begin{equation}
\label{equ:control:system}
\begin{tikzcd}
\mathcal{U} \arrow[r, "f"] \arrow[rd, "\pi_u"] & TX \arrow[d, "\pi_x"]                       \\
                                  & X \arrow[lu, "\sigma", dotted, bend left=49]
\end{tikzcd}.
\end{equation}
Representing a control system through \eqref{equ:control:system} reveals how state-dependent input constraints come to play when synthesizing (stabilizing) feedback controllers. Studying control systems through such a diagram is commonly attributed to Brockett (\textit{e.g.}, see \cite[\S 13.5]{ref:nonlinear1990NvdS} and the accompanying notes). Note that \eqref{equ:control:system} generalizes the case of $\dot{x}=f(x,u)$ with $f:X\times U\to TX$, that is, for such a control system we work with a trivial vector bundle $\pi_u:\mathcal{U}=X\times U\to X$. This remark is important as one might be accustomed to the notation $f(x,u)$ and not $f(u)$ (or $f(u_x)$). Also note that one typically\footnote{Global sections do not always exist. For example, consider what happens when one removes the zero section from the tangent bundle of the $2$-sphere, that is, $T\mathbb{S}^2\setminus 0_{T\mathbb{S}^2}\to \mathbb{S}^2$; a global section would contradict the well-known ``\textit{hairy ball}" theorem.} studies \textit{local} sections, that is, a map $\sigma:N\subseteq X\to \mathcal{U}$ such that $\pi_u\circ \sigma = \mathrm{id}_N$.

Next, we provide both (classes of) necessary conditions, for a particular feedback stabilization problem. Here, $0_{TX}:X\to TX$ denotes the \textit{zero section} of the tangent bundle $TX$, that is, $x \mapsto 0_{TX}(x)=(x,0)$, locally. With some abuse, we will also call the image $0_{TX}(X)$ the zero section. 

\subsection{Necessary conditions for feedback stabilization}
\label{sec:nec:cond:1.2}
Suppose we work with a control system $\mathcal{C}=(X,\mathcal{U},f)$. Then, given a sufficiently regular (\textit{e.g.}, such that an Euler characteristic\footnote{The Euler characteristic (\textit{e.g.}, see \cite[p. 116]{ref:guillemin2010differential}) is a topological invariant that constrains, in particular, dynamical behaviour. For instance, let $M$ be a smooth, compact manifold without boundary, then $\chi(M)=0$ if and only if $M$ admits a non-vanishing vector field \cite[p. 146]{ref:guillemin2010differential}, \textit{e.g.}, $\chi(\mathbb{S}^1)=0$ while $\chi(\mathbb{S}^2)=2$. See~\eqref{equ:euler:charac:homology} for the definition we use.} $\chi(\cdot)$ is well-defined) compact subset $A\subseteq X$ one may ask if there is a (local) section $\sigma$ such that $A$ is (locally) asymptotically stable under the closed-loop vector field $f\circ \sigma$. This is a question of stabilization through \textit{continuous feedback} (see \S\ref{sec:Lyapunov stability} for more on stability theory). A necessary condition that generalizes Brockett's condition \cite[Thm. 1.(iii)]{ref:brockett1983asymptotic}, as proposed by Kvalheim and Koditschek \cite[Thm. 3.2]{ref:kvalheim2021necessary}, is as follows. Note, one needs $\chi(A)\neq 0$.  
\\\\
\textbf{\mylabel{adversary:condition}{(i)} The adversary condition} For any neighbourhood $N\subseteq X$ of $A$, there is a neighbourhood $N_0\subseteq TX$ of the zero section $0_{TX}$, such that for any $N_0$-valued continuous vector field $E$ over $N$, $E(x) = f(u_x)$ for some $x\in N$ and $u_x\in \pi_u^{-1}(N)$.
\\\\
Observe that, like Brockett's original necessary condition, the adversary condition hinges on continuity, that is, around the attractor $A$, the closed-loop vector field can be made arbitrarily small in magnitude. As we argue below, this is not necessary when one is concerned with stability and continuity of solutions (\textit{e.g.}, semiflows) alone; continuity of the closed-loop vector field is not necessary.

Next, Coron's condition has been recently generalized by Kvalheim \cite[Thm. 2]{ref:kvalheim2022obstructions}. Essentially, by generalizing the homotopy argument as exploited by Krasnosel'ski\u{\i} and Zabreiko \cite[Lem. 52.1]{ref:krasnosel1984geometrical}, \textit{cf}. \cite[Thm. 1]{ref:kvalheim2022obstructions}. Regarding notation, given a subset $Y\subseteq X$, let $\Sigma_Y:= \pi_u^{-1}(Y)\setminus f^{-1}(0_{TX})$, that is, $\Sigma_Y\subseteq \mathcal{U}$ is the subset of the control space not mapped to the zero section of $T_Y X:=\cup_{y\in Y}T_y X$ under $f$, \textit{cf}. \cite{ref:Coron1990}.  
\\\\
\textbf{\mylabel{homology:condition}{(ii)} The homology condition}
\label{(ii)} Let $G$ be any uniquely integrable continuous vector field on $X$ such that $A$ is asymptotically stable under $G$. Then, for any sufficiently small neighbourhood $N$ of $A$ we have
\begin{align}
\label{equ:homology:cond}
    G_* H_{\bullet} (N\setminus A) \subseteq f_* H_{\bullet}(\Sigma_{N\setminus A}) \quad \text{in}\quad  H_{\bullet} (T_{N\setminus A}X\setminus 0_{TX}).
\end{align}
Here, $H_{\bullet}(X)$ denotes graded homology (\textit{i.e.}, $H_{\bullet}(X)=\oplus_k H_k(X)$). In this note, we employ singular homology with integer coefficients for simplicity of exposition. One may consider \v{C}ech cohomology instead to capture \textit{shape} equivalence of attractors and their basin of attraction. For more details regarding this generalization, see \cite{ref:kapitanski2000shape,ref:gobbino2001topological,ref:kvalheim2021necessary}.  

When it comes to comparing these conditions, Kvalheim shows the following. When the control system has the structure of a vector bundle, that is, $\pi_u:\mathcal{U}\to X$ is a vector bundle, and $X$ is orientable, then, the homology condition is a sharper necessary condition (\textit{i.e.}, {\myref{homology:condition}{(ii)}}${\,\,\Rightarrow}$ {\myref{adversary:condition}{(i)}}: when the homological condition holds, so does the adversarial condition) \cite[Thm. 4]{ref:kvalheim2023relationships}. 

Before continuing, we recall some terminology. Let $E$ and $B$ be topological spaces, then, a tuple $(E,B,\pi)$, such that $\pi:E\to B$ is a continuous surjective map, with a vector space structure on each fiber $\pi^{-1}(b)$, is a $k$-dimensional \textbf{\textit{vector bundle}} when for any $b\in B$, there is an open neighbourhood $U\ni b$ such that $\pi^{-1}(U)$ is homeomorphic to $U\times \mathbb{R}^k$, with this homeomorphism  restricting to a vector space isomorphism when taking $\pi^{-1}(b)$ to $\{b\}\times \mathbb{R}^k$. We point the reader to \cite{ref:HatcherVBK2003} for more on vector bundles. Relaxing the vector space structure provides us with a fiber bundle. That is, let $E$, $B$ and $F$ be topological spaces, then the tuple $(F,E,B,\pi)$ is a \textbf{\textit{fiber bundle}} when $\pi:E\to B$ is a continuous surjective map and for any $b\in B$ there is an open neighbourhood $U\ni b$ such that $\pi^{-1}(U)$ is homeomorphic to $U\times F$, and this homeomorphism is fiber-preserving.  

Now, let $\pi:E\to B$ be an $n$-dimensional vector bundle (\textit{e.g.}, $\pi_x:TB\to B$), then, there are several methods to define the so-called \textbf{\textit{sphere bundle}} $\mathbb{S}^{n-1}\hookrightarrow E \to B$.
One may simply demand that all fibers are topological spheres, that is, $\pi^{-1}(b)\cong \mathbb{S}^{n-1}$ for all $b\in B$. Indeed, a sphere bundle is a fiber bundle that is not a vector bundle. As vector bundles, with $B$ being paracompact and Hausdorff, have an inner product \cite[Prop. 1.2]{ref:HatcherVBK2003}, one may also define the sphere bundle, induced by $E$, as $S(E):=\{e\in E:\langle e,e\rangle =1,e\in \pi^{-1}(b),b\in B\}$. We will use the latter. Note that if $B$ is a smooth manifold with tangent bundle $TB$, we may endow $B$ with a Riemannian metric $g$, and thus, $S(TB)$ will be defined using this metric. We will only comment on the employed metric if relevant. The norm induced by this metric will be denoted by $\|\cdot \|_g$ or simply $\|\cdot \|$. In general, one could define the zero section $0_E:B\to E$ and study a quotient space of $E\setminus 0_E$, \textit{e.g.}, use $\mathbb{R}^n\setminus \{0\}\cong \mathbb{S}^{n-1}$. For instance, regarding \eqref{equ:homology:cond} we have that
\begin{equation}
\label{equ:tangent:sphere:bundle}
    T_{N\setminus A}X\setminus 0_{TX} \cong S(T_{N\setminus A}X).
\end{equation}

Then, as stated before, up and including (oriented) vector bundles, it is known that the homology condition implies the adversary condition. One may wonder if the comparison stops there.  
To that end, in the case of fiber bundles, Kvalheim provides the following counterexample \cite[Ex.4]{ref:kvalheim2023relationships}. Counter in the sense that this example shows that the homology condition~{\myref{homology:condition}{(ii)}} does not always imply the adversary condition~{\myref{adversary:condition}{(i)}} when looking at fiber bundles that are not vector bundles. 

\begin{example}[{\myref{homology:condition}{(ii)}}${\,\,\nRightarrow}$ {\myref{adversary:condition}{(i)}}: a counterexample using sphere bundles {\cite[Ex.4]{ref:kvalheim2023relationships}}]
\label{ex:counter}
\upshape{
    Let $X:=\mathbb{R}^n$, $\pi_u:\mathcal{U}\to X$ be the unit-sphere bundle, that is, $\mathcal{U}:=S(TX)\subseteq TX$, and let $f:\mathcal{U}\hookrightarrow TX$ be the inclusion map. Now, consider $A:=\{0\}$. The homology condition~{\myref{homology:condition}{(ii)}} readily holds. To see this, first note that (any) $G$ is homotopy equivalent to $G/\|G\|$, away from $A$. Then pick $\sigma:=G/\|G\|$ and as $G=f\circ \sigma$ conclude on set inclusion. The adversary condition~{\myref{adversary:condition}{(i)}}, on the other hand, fails as there is no input available to render the vector field sufficiently small, that is, for any vector field $E$ such that $\|E(x)\|\in [0,1)$ for all $x\in X$, there is no admissible $u_x\in \pi_u^{-1}(U)$ as $\|f(u_x)\|=1\notin (0,1)$.
    }
    \exampleEnd
\end{example}

Now consider the following example.
\begin{example}[Finite-time stability and sphere bundles]
\label{ex:finite:time:sphere}
\upshape{
Consider a control system on $\mathbb{R}^n$, where $\dot{x}=u$ with $u\in \mathbb{S}^{n-1}$, that is, we are in the setting as in Example~\ref{ex:counter} above. Now define the feedback $\mu$ as follows. On $\mathbb{R}^n\setminus \{0\}$ 
\begin{align*}
    x\mapsto \mu(x) := -\frac{x}{\|x\|_2}, 
\end{align*}
whereas at $x=0$, we let $x\mapsto \mu(x)= \mathbb{S}^{n-1}$, that is, the map is set-valued. 
Now, we can study $\dot{x}\in \mathcal{F}[\mu](x)$, for $\mathcal{F}$ the Filippov operator (see \S\ref{sec:Lyapunov stability}). Indeed, the origin is (finite-time) asymptotically stable under this dynamical system\footnote{The corresponding semiflow is given by $\mathbb{R}_{\geq 0}\times \mathbb{R}^n\ni (t,x_0)\mapsto \varphi^t(x_0):= (1-t/\|x_0\|_2)x_0$ for $t\leq \| x_0\|_2$ and $\varphi^t(x_0):=0$ otherwise. Note, this map is continuous.}, despite failing to satisfy the adversary condition~{\myref{adversary:condition}{(i)}}.}
\exampleEnd
\end{example}

The example above shows that the counterexample (Ex.~\ref{ex:counter}) only holds due to a lack of continuity of any closed-loop vector field around $0$. Stability of the origin under these vector fields, on the other hand, is not obstructed (see \S\ref{sec:discussion}(c) for a similar example in the context of CBFs). Neither is continuity of the solutions. Still, is the fact that the homology condition \textit{does} hold for Example~\ref{ex:counter} a coincidence or a feature? Better yet, can we alter the adversarial condition and make it suit Example~\ref{ex:counter}? With this in mind, we ask the following.
\\\\
\textbf{Question B:} ``\textit{Is there an adversarial condition, suitable for Example~\ref{ex:counter}-\ref{ex:finite:time:sphere}; a condition that is applicable in a sufficiently general setting; a condition that is again implied by the homology condition?}"
\\\\
An affirmative answer would resolve Question A as motivated and stated in the above---the homology condition should be understood as being simply more general than the adversarial condition. Indeed, it turns out that a mild adjustment of the adversary condition solves this question, as formalized in this note. 
It is clear that for such a generalization, the \textit{magnitude} of adversaries should not play an important role, their ``\textit{direction}" is of importance. Indeed, in Example~\ref{ex:finite:time:sphere}, the control system is such that for any neighbourhood $N$ of the point $0$, one can find for any $z\in \mathbb{S}^{n-1}$ some $u_x\in \pi_u^{-1}(N\setminus \{0\})$ such that $z=f(u_x)/\| f(u_x)\|_2$.

\subsection{Control systems revisited}
\label{sec:discont:control:system}
With Example~\ref{ex:finite:time:sphere} in mind, we introduce a slight modification into the definition of a control system $\mathcal{C}=(X,\mathcal{U},f)$. Let $D\subseteq X$ be a set where closed-loop vector fields may be discontinuous, but remain measurable and locally essentially bounded. If $D\neq \varnothing$, we interpret the resulting discontinuous differential equation in the sense of Filippov, that is, we study set-valued maps $\mathcal{F}[f\circ \sigma]:X\rightrightarrows TX$ (see \S\ref{sec:Lyapunov stability}). Note, the discontinuity may arise due to $\sigma$, $f$ or both. For simplicity, however, we will assume that $f$ and $\sigma$ are only allowed to be discontinuous on $\Sigma_D$ and $D$, respectively (and thus, those maps are continuous on $\Sigma_{X\setminus D}$ and $X\setminus D$, respectively). Overloading notation, we denote the resulting tuple by $\mathcal{C}=(X,\mathcal{U},f,D,\mathcal{F})$ (\textit{i.e.}, if $D=\varnothing$, we recover $(X,\mathcal{U},f)$ as in that case $\mathcal{F}[f\circ \sigma]=f\circ \sigma$).   

We remark that this definition is inspired by Example~\ref{ex:counter}-\ref{ex:finite:time:sphere} and finite-time stabilization. One may study a more general setting. 

\section{Lyapunov stability}
\label{sec:Lyapunov stability}
In this section we follow Mayhew and Teel~\cite{ref:mayhew2011topological} and highlight a converse Lyapunov theorem applicable to some discontinuous dynamical systems, see also \cite{ref:clarke1998asymptotic}. A reader familiar with stability theory may safely skip this section.

Let $X$ be an $n$-dimensional smooth, Hausdorff and second-countable manifold (see \cite[Ch. 1]{Lee2}). Then, there is a smooth embedding $X\hookrightarrow \mathbb{R}^{2n+1}$ ($X$ is diffeomorphic to some closed submanifold of $\mathbb{R}^{2n+1}$ \cite[Cor. 6.16]{Lee2}). Hence, we will assume in what follows that $X$ is an embedded submanifold of some Euclidean space.

Now, we will study stability under set-valued vector fields $F:X\rightrightarrows TX$, that is, $F(x)\subseteq T_x X$. More concretely, we study stability under the differential inclusion 
\begin{equation}
    \label{equ:DI}
    \dot{x}\in F(x),\quad x\in X. 
\end{equation}
To handle the inclusion \eqref{equ:DI}, it is convenient to assume that the map $F:X\rightrightarrows TX$ satisfies the so-called \textbf{\textit{basic conditions}} \cite[Def. 5]{ref:mayhew2011topological}. These conditions are satisfied precisely when
\begin{enumerate}[(i)]
    \item $F$ is \textit{locally bounded}, that is, for any compact set $K_1\subseteq X$ there is a compact set $K_2\subseteq TX$ such that $F(K_1)\subseteq K_2$;
    \item $F$ is \textit{outer semicontinuous}, that is, $\mathrm{graph}(F)$ is closed; and
    \item $F(x)$ is nonempty, compact and convex, for any $x\in X$. 
\end{enumerate}
Recall, if $F:X\to TX$ is not set-valued, but still locally bounded, then $\mathrm{graph}(F)$ being closed implies that $F$ is continuous. 

Consider now a possibly discontinuous ODE $\dot{x}=f(x)$, $x\in X$, with $f$ measurable and locally essentially bounded. Then, the \textbf{\textit{Filippov operator}} $\mathcal{F}$ \cite[p. 85]{ref:filippov1988differential} is defined pointwise through
\begin{align}
\label{equ:Filippov}
    \mathcal{F}[f](x) := \bigcap_{r>0} \bigcap_{N:\mu(N)=0} \overline{\mathrm{conv}}f (\mathbb{B}_r(x)\setminus N),
\end{align}
where $\mu$ denotes the appropriate Lebesgue measure and $\overline{\mathrm{conv}}(B)$ the closure of the convex hull of a set $B$. Applying the Filippov operator to such an $f$ leads to a set-valued vector field that satisfies the basic conditions. The rationale is to pass from a discontinuous differential equation to a more regular differential inclusion. 

Now, a solution to an inclusion \eqref{equ:DI} is an absolutely continuous curve $[0,T]\ni t\mapsto \xi(t)\in X$ that satisfies
\begin{align*}
    \left.\frac{\mathrm{d}}{\mathrm{d}t}\xi(t)\right|_{t=s} \in F(\xi(s)) \quad \text{for a.e.}\, s\in [0,T]. 
\end{align*}
Then, an inclusion \eqref{equ:DI} is said to be \textit{forward complete} when all solutions are defined for all $t\geq 0$. 

Regarding stability, a compact set $A\subseteq X$ is said to be (locally) \textbf{\textit{asymptotically stable}} under (a forward complete) inclusion \eqref{equ:DI} when
\begin{enumerate}[(i)]
     \item $A$ is \textit{stable}, that is, for any open neighbourhood $U_{\infty}$ of $A$, there is an open neighbourhood $U_0$ of A such that solutions $t\mapsto \xi(t)$ under \eqref{equ:DI} that start in $U_0$ stay in $U_{\infty}$ for all time (\textit{i.e.}, $\xi(0)\in U_0$ $\implies$ $\xi(t)\in U_{\infty}\, \forall t\geq 0$); and
    \item $A$ is \textit{attractive}, that is, there is an open neighbourhood $V$ of $A$ such that for any open neighbourhood $W$ of $A$ there is $T\geq 0$ such that solutions $t\mapsto \xi(t)$ under \eqref{equ:DI} starting in $V$ enter $W$ in time $T$, and stay there (\textit{i.e.}, $\xi(0)\in V$ $\implies \xi(t)\in W\,\forall t\geq T$).  
\end{enumerate}
Let $A$ be asymptotically stable, then the \textit{basin of attraction}, denoted $B(A)$, is the union of all neighbourhoods $V$, such that $A$ is attractive from $V$. As $A$ is compact, without loss of generality, one may define stability through metric neighbourhoods as well \cite[Ch. V]{ref:bhatia1970stability}.  

Then, based on \cite{ref:cai2008smooth2}, \cite{ref:mayhew2011topological} provides the following converse Lyapunov theorem. 

\begin{theorem}[A converse Lyapunov theorem for inclusions {\cite[Cor. 13]{ref:mayhew2011topological}}]
    \label{thm:converse:Lyapunov:DI}
    Let the compact subset $A$, of some embedded submanifold $X\subseteq \mathbb{R}^m$, be asymptotically stable under $\dot{x}\in F(x)$, with $F:X\rightrightarrows TX$ satisfying the basic conditions, then, there is a smooth function $V:B(A)\to \mathbb{R}_{\geq 0}$ such that
\begin{enumerate}[(i)]
    \item $\alpha_1(\omega(x))\leq V(x)\leq \alpha_2(\omega(x))$, where $\alpha_1,\alpha_2\in \mathcal{K}_ {\infty}$ and $\omega$ is a any proper continuous indicator function for $A$; and
    \item $\max_{v\in F(x)} \langle \nabla V(x),v\rangle \leq - V(x)$.
\end{enumerate}
\end{theorem}

In the next section we exploit Item (ii) from Theorem~\ref{thm:converse:Lyapunov:DI}. It allows for capturing the ``\textit{direction}'' of a discontinuous vector field via properties of a smooth object (a gradient flow generated by $V$).

\textbf{Homotopy to a smooth gradient flow} Let $F$ be a vector field on a smooth manifold $X$ such that $A$ is asymptotically stable with basin of attraction $B(A)\subseteq X$. Let $F$ be continuous on $B(A)\setminus A$ and suppose there is a smooth Lyapunov function $V$ such that $L_F V< 0$ on $B(A)\setminus A$ (here, $L_F V$ denotes the Lie derivative of $V$ along $F$ \cite[p. 227]{Lee2}). Since $X$ is smooth, $X$ admits a Riemannian metric $g$ \cite[Prop. 3.13]{Lee2}. Then, we can define the \textbf{\textit{Riemannian gradient}} of $V$ through the relation $\mathrm{D}f(x)[v]=\langle \mathrm{grad}\, V(x),v\rangle_g$; this relation must hold for any smooth function $f:B(A)\to \mathbb{R}$ and any $v\in T_x X$. It follows that the set $A$ is also asymptotically stable under the vector field $-\mathrm{grad}\,V$ (\textit{i.e.}, check that $ -\mathrm{grad}\, V(V)<0$ on $B(A)\setminus A$). Now, define for any $s\in [0,1]$ the vector field $H_s := (1-s)F-s\mathrm{grad}\,V$. It follows that $H_s$ is continuous on $B(A)\setminus A$ and that $L_{H_s}V<0$ on $B(A)\setminus A$, as such, $[0,1] \ni s\mapsto H_s$ yields a homotopy, on $B(A)\setminus A$, from $F$ to a smooth gradient flow. In particular, it follows that $H_s$ is nonzero for any $s$, that is, we work with the continuous map
\begin{equation*}
    H_s: B(A)\setminus A \to T_{B(A)\setminus A} X \setminus 0_{TX} 
\end{equation*}

As we will work in a setting where the existence of such a Lyapunov function can be guaranteed (\textit{e.g.}, via \cite{ref:mayhew2011topological,ref:fathi2019smoothing}), we can and will use the existence of this homotopy in the remainder of this note. For simplicity, we will write $\nabla V$ instead of $\mathrm{grad}\, V$. 

\section{A spherical adversary condition}
\label{sec:main}
The first result is a Brockett-like necessary condition for stabilization. It will allow for correctly interpreting Examples~\ref{ex:counter} and \ref{ex:finite:time:sphere} in relation to the homology condition. 

The result is inspired by \cite{ref:kvalheim2021necessary}, but with the focus on a less regular closed-loop vector field. This is possible by simply focussing on the smooth gradient flow generated by a Lyapunov function. For simplicity, our compact attractors are slightly less general than the ones in \cite[\S{2}]{ref:kvalheim2021necessary} (\textit{e.g.}, Kvalheim and Koditschek consider \textit{shape} equivalence of an attractor $A$ and its basin of attraction of $B(A)$, we consider \textit{homotopy} equivalence), as captured by the following assumption.


\begin{assumption}[On the pair $(A,X)$, with $A\subseteq X$]
\label{ass:AX}
Let $X$ be a smooth, orientable, second-countable, Hausdorff manifold and let $A$ be a compact subset of $X$ such that $\iota:A\hookrightarrow X$ is cofibration (\textit{e.g.}, $A$ is an embedded submanifold). 
\end{assumption}
Note that the pair $(A,X)$ from Example~\ref{ex:finite:time:sphere} satisfies Assumption~\ref{ass:AX}. We point the reader to \cite[Ch. VII]{ref:bredon1993topology} for more comments on cofibrations, and to \cite{ref:jongeneelECC24,ref:jongeneel2025topological} for more on the relation with dynamical systems. Now, given a set $A$, we define its \textbf{\textit{Euler characteristic}} through (see \cite[Thm. 2.44]{Hatcher}): 
\begin{align}
\label{equ:euler:charac:homology}
    \chi(A) := \sum^{\infty}_{k=0}(-1)^k\,\mathrm{rank}\,H_k(A).
\end{align}
When the pair $(A,X)$ satisfies Assumption \ref{ass:AX} and $A$ is asymptotically stable under some sufficiently regular differential equation, we know that \eqref{equ:euler:charac:homology} is well-defined, \textit{e.g.}, since $A$ is homotopy equivalent to a smooth manifold, see \cite{ref:jongeneelECC24}.

As will be clarified, it is useful to \textit{momentarily} impose another, stronger assumption. This assumption is akin to \cite{ref:mansouri2007local,ref:mansouri2010topological} and helps in making sense of ``\textit{directions}".
\begin{assumption}[On $X$]
\label{ass:X}
Let $X$ be an open subset of $\mathbb{R}^n$. 
\end{assumption}
The most important ramification of Assumption~\ref{ass:X} is that $TX\cong X \times \mathbb{R}^n$. Hence, (closed-loop) vector fields $F:X\to TX$, or other fiber preserving maps into $TX$, can be identified with maps $F':X\to \mathbb{R}^n$. This is not the most general setting we will consider; we comment on appropriate generalizations later (\textit{e.g.}, see Corollary~\ref{cor:homology} and \S\ref{sec:generalizaions}).

\begin{theorem}[A spherical adversary condition]
\label{thm:spherical:adversary}
    Let $(A,X)$ satisfy Assumptions~\ref{ass:AX}-\ref{ass:X} and consider a control system $\mathcal{C}=(X,\mathcal{U},f,A,\mathcal{F})$ in the sense of \S\ref{sec:discont:control:system}. Let $N\subseteq X$ be an open neighbourhood of $A$ and suppose there is a section $\sigma:N\to \mathcal{U}$ such that $A$ is asymptotically stable under $f\circ \sigma$, in the sense of Filippov. 
      Then,  
      \begin{enumerate}[(i)]
          \item $\chi(A)$ is well-defined; and
          \item  if $\chi(A)\neq 0$, then, for any open neighbourhood $N'\subseteq N$ of $A$ 
      \begin{align*}
          \Sigma_{N'\setminus A}\ni u_x \mapsto \frac{f(u_x)}{\|f(u_x)\|_2}\in \mathbb{S}^{n-1}
      \end{align*}
      is surjective, with $f$ understood as $f:\mathcal{U}\to \mathbb{R}^n$. 
      \end{enumerate}
\end{theorem}
\begin{proof}
By Assumption~\ref{ass:AX}, we may assume that $X$ is an embedded manifold of some Euclidean space. 

Then, as $A$ is asymptotically stable under a closed-loop vector field $\mathcal{F}[f\circ \sigma]:N\subseteq X\rightrightarrows T_N X$ that satisfies the basic conditions (by the definition of $\mathcal{C}$), Theorem~\ref{thm:converse:Lyapunov:DI} (see \cite[Cor. 13]{ref:mayhew2011topological}) implies there exists a $C^{\infty}$ Lyapunov function $V:B(A)\subseteq N \to \mathbb{R}_{\geq 0}$. In particular, $A$ is asymptotically stable under $\dot{x}=-\nabla V(x)$. Let $M_c:= V^{-1}([0,c])\subseteq B(A)$ for some $c>0$. 

As $X$ is orientable, $A$ is compact and $\nabla V$ is of constant rank on $\partial M_c$, the set $\partial M_c$ is a compact, orientable, $(n-1)$-dimensional embedded submanifold of $X$ \cite[Thm. 5.12]{Lee2}. It may happen that $\partial M_c$ is not connected (\textit{e.g.}, let $X=\mathbb{R}^3$ and $A=\mathbb{S}^2$, see also \cite{ref:mansouri2007local}). However, by compactness, $\partial M_c$ consists out of finite components, say $k\in \mathbb{N}_{> 0}$ components. 

Then, under our standing assumptions (\textit{i.e.}, $A$ is asymptotically stable and $\iota:A\hookrightarrow X$ is a cofibration), we have that $\chi(M_c)$ is well-defined\footnote{Construct the closed manifold $2M_c$ via glueing a copy of $M_c$ to itself along $\partial M_c$. Hence, $\chi( 2 M_c) = 2\chi(M_c)-\chi(\partial M_c)$. Then, if $\partial M_c$ is odd-dimensional, we have that $\chi(M_c)=\frac{1}{2}\chi(2M_c)$, otherwise, $\chi(M_c)=\frac{1}{2}\chi(\partial M_c)$.} and $\chi(M_c)=:\chi(A)$, as $B(A)$, $A$ and $M_c$ are all homotopy equivalent \cite{ref:jongeneelECC24}. 


As $TX\cong X\times \mathbb{R}^n$, one can identify sections into the sphere bundle $S(TX)$ with maps into $\mathbb{S}^{n-1}$. Now consider the maps $g,h:\partial M_c\to \mathbb{S}^{n-1}$ as defined through $x\mapsto g(x):=  f(\sigma(x))/\| f(\sigma(x))\|_2$ and $x\mapsto h(x):= -\nabla V(x)/\| \nabla V(x)\|_2$. Note that $h$ is smooth, but $g$ is only known to be continuous. Then, as $L_{f\circ \sigma} V<0$ on $x\in B(A)\setminus A$ (and thus on $\partial M_c$), we can construct a (single) straight-line (continuous) homotopy from $g$ to $h$ (recall the last paragraph of \S\ref{sec:Lyapunov stability}). This means that $\mathrm{deg}(h)=\mathrm{deg}(g)$~\cite[p. 108]{ref:guillemin2010differential}. 
    Then, as $\chi(A)\neq 0$, the Poincar\'e-Hopf theorem \cite[p. 35-36]{ref:milnor65} implies
    \begin{align*}
        (-1)^n\mathrm{deg}(h)=\mathrm{deg}(-h) =\chi(A)\neq 0,
    \end{align*}
and thus we know that $\mathrm{deg}(g)=\mathrm{deg}(h)\neq 0$. Hence, $g$ must be surjective over $\partial M_c$, where our choice of $c>0$ was arbitrary. Then, to conclude, given any open neighbourhood $N'\subseteq N$ of $A$, we can always find a sufficiently small $c'>0$ such that $M_{c'}\subseteq N'$ (\textit{e.g.}, as $X$ is smooth, endow $X$ with a Riemannian metric and appeal to standard results regarding metric neighbourhoods of compact sets \cite[\S27]{ref:munkrestopology}). 
\end{proof}

As discussed in \S\ref{sec:discussion}(b), $\chi(A)\neq 0$ is not an artifact of the proof, this condition cannot be removed. We also emphasize that Theorem~\ref{thm:spherical:adversary} is not equivalent to showing that the Gauss map from some smooth levelset is surjective, as pointed out in \cite[Footnote 4]{ref:kvalheim2021necessary}, this is always the case. Instead, we focus on the degree. The crux is that surjectivity is not homotopy invariant, but the degree of a map is.

Now, given Theorem~\ref{thm:spherical:adversary}, we study the following version of the adversary condition.  
\\\\
\textbf{\mylabel{spherical:adversary:condition}{(iii)} The spherical adversary condition} For any open neighbourhood $N\subseteq X$ of $A$ and any $z\in \mathbb{S}^{n-1}$, there is a $u_x\in \pi_u^{-1}(N\setminus A)$ such that $z = f(u_x)/\|f(u_x)\|_2$, with $f$ understood as $f:\mathcal{U}\to \mathbb{R}^n$.
\\\\
Indeed, the spherical adversary condition applies to Examples~\ref{ex:counter}-\ref{ex:finite:time:sphere} whereas the standard adversary condition does not. 

To continue, we show that the spherical adversary condition complies with what is known regarding the nonholonomic integrator, \textit{e.g.}, see \cite{ref:brockett1983asymptotic,ref:zabczyk1989,ref:ryan1994brockett}. 
\begin{example}[The nonholonomic integrator]
\label{ex:nonholonomic}
\upshape{
    Consider the classical nonholonomic integrator
    \begin{equation}
        \label{equ:nonholonomic}
         \begin{aligned}
        \dot{x}_1 & = u_1\\
        \dot{x_2} &= u_2\\
        \dot{x}_3 &= x_1 u_2 - x_2 u_1,
    \end{aligned}
    \end{equation}
    compactly written as $\dot{x}=f(x,u)=:\sum_i g_i(x) u_i$. Interestingly, this system is controllable (\textit{e.g.}, see that \eqref{equ:nonholonomic} is driftless and check that $\mathrm{rank}(g_1,g_2,[g_1,g_2])=3$ \cite[Ch. 3]{ref:nonlinear1990NvdS}). Now consider stabilizing the origin, it readily follows that Brockett's condition fails as there is no pair $(x_{\varepsilon},u_{\varepsilon})$ such that $f(x_{\varepsilon},u_{\varepsilon})=(0,0,\varepsilon)$, for some $\varepsilon>0$. Hence, there is no differentiable (or even continuous) feedback such that the origin is stabilized. Regarding Coron's condition, one may check that $f_* H_{n-1}(\Sigma) \cong f_* H_2(\mathbb{S}^1) \not\cong \mathbb{Z}$. Clearly, there is also no pair $(x_s,u_s)$ such that $f(x_s,u_s)/\| f(x_s,u_s)\|_2=(0,0,1)$. Thus, the spherical adversary condition fails to hold, that is, introducing mild discontinuities does not resolve the question of stabilization. However, we may stabilize a set $A$, in finite-time, with $\chi(A)=0$, as shown in \S\ref{sec:discussion}(b).   
    } \exampleEnd
\end{example}

The spherical adversary condition {\myref{spherical:adversary:condition}{(iii)}} is all about the ``\textit{directions}" a control system can generate around some set $A$, and the intention of proposing precisely this condition is that it is directly of use (\textit{e.g.}, there is no need to first perform any form of regularization via Filippov operators and the like). 

A sidenote here is that in some cases, an application of the Filippov operator fills in precisely the gap in the tangent space that would otherwise prevent the standard adversary condition from being true (\textit{e.g.}, going from a sphere- to disk bundle).

\begin{remark}[Comparison to~\cite{ref:ryan1994brockett}]
    \upshape{Ryan studies a setting akin to Brockett (\textit{i.e.}, for $\dot{x}=f(x,u)$, with $f$ continuous, find a feedback to stabilize $0$), but Ryan allows for mildly discontinuous feedback. Here, the additional assumption is that for $K$ convex, $f(x,K)$ is convex. Moreover, feedback controllers are upper semicontinuous maps $x\mapsto \mu(x)$, such that $\mu(x)$ is nonempty, compact and convex, with $0\in \mu(0)$. A Brockett-like necessary condition can be recovered \cite[Cor. 1]{ref:ryan1994brockett}. This construction is reminiscent of applying the Filippov operator to a discontinuous vector field, however, it is more stringent (\textit{i.e.}, compare applying the Filippov operator to either a discontinuous vector field or only to the controller, see \cite{ref:ceragioli2002some}). A simple example that cannot be handled using the machinery from \cite{ref:ryan1994brockett} would be $\dot{x}=\sin(u)$ with $u\in \mathbb{S}^0$. Nonetheless, $\dot{x}=\sin(-\mathrm{sgn}(x))$ renders the origin asymptotically stable, in the sense of Filippov.}
    \exampleEnd
\end{remark}

Next, before comparing the spherical adversary condition {\myref{spherical:adversary:condition}{(iii)}} to the homology condition {\myref{homology:condition}{(ii)}}, we ensure that the homology condition is still meaningful (\textit{i.e.}, still necessary) in our setting, that is, in the context of Section~\ref{sec:discont:control:system}. This readily follows from Kvalheim's homotopy theorem \cite[Thm. 2]{ref:kvalheim2023relationships}, \cite[Thm. 1]{ref:kvalheim2022obstructions} and the structure as exploited in the first part of the proof of Theorem~\ref{thm:spherical:adversary}. Note, we do not rely on Assumption~\ref{ass:X} or $\chi(A)\neq 0$.

\begin{corollary}[The homology condition under discontinuous feedback]
\label{cor:homology}
    Let $(A,X)$ satisfy Assumption~\ref{ass:AX} and consider a control system $\mathcal{C}=(X,\mathcal{U},f,A,\mathcal{F})$ in the sense of \S\ref{sec:discont:control:system}. Let $N\subseteq X$ be an open neighbourhood of $A$ and suppose there is a section $\sigma:N\to \mathcal{U}$ such that $A$ is asymptotically stable under $f\circ \sigma$, in the sense of Filippov, then the homology condition holds. 
\end{corollary}
\begin{proof}
Let $A$ be asymptotically stable on some sufficiently small neighbourhood $N'\subseteq N$ of $A$. First, one constructs a (straight-line) homotopy from $f\circ \sigma$ to $-\nabla V$, without introducing any zeros, on $N'\setminus A$. Then, by Kvalheim's homotopy theorem \cite[Thm. 1]{ref:kvalheim2022obstructions}, there is a homotopy from $-\nabla V$ to $G$, on $W\setminus A$, where $W\subseteq N'$ is a possibly smaller neighbourhood of $A$, again, without introducing any zeros. Now we may conclude as in the proof of \cite[Thm. 3]{ref:kvalheim2023relationships}, that is
\begin{align*}
    G_*H_{\bullet}(W\setminus A)=-\nabla V_* H_{\bullet}(W\setminus A) = (f\circ \sigma)_* H_{\bullet}(W\setminus A)\subseteq f_* H_{\bullet}(\Sigma_{W\setminus A}),
\end{align*}
in $H_{\bullet}(T_{W\setminus A} \setminus 0_{TX})$. Note, here we do exploit that both $f$ and $\sigma$ are continuous on $\Sigma_{W\setminus A}$ and $W\setminus A$, respectively. For the spherical adversary condition this was not needed (only continuity of $f\circ \sigma$ on $W\setminus A$ was used, and only indirectly to appeal to smooth Lyapunov functions and the topological degree). 
\end{proof}

\begin{theorem}[{\myref{homology:condition}{(ii)}}${\,\,\Rightarrow}$ {\myref{spherical:adversary:condition}{(iii)}}: the homology condition is always stronger than the spherical adversary condition]
\label{thm:compare}
Let the pair $(A,X)$ satisfy Assumptions~\ref{ass:AX}-\ref{ass:X} and consider a control system $\mathcal{C}=(X,\mathcal{U},f,A,\mathcal{F})$ in the sense of \S\ref{sec:discont:control:system} with $\chi(A)\neq 0$. Then, if the homology condition holds, so does the spherical adversary condition.  
\end{theorem}
\begin{proof}
As in \cite{ref:kvalheim2023relationships}, we show that if the spherical adversary fails to hold, so does the homology condition. Note, we use again the identification $TX\cong X\times \mathbb{R}^n$, \textit{e.g.}, $H_{\bullet}(W\times \mathbb{R}^n)\cong H_{\bullet}(W)$.

    With the above in mind, consider $G$ as in the homology condition. By the properties of $G$, we know there is a smooth Lyapunov function $V$ (with respect to $A$ being asymptotically stable under $G$) \cite{ref:fathi2019smoothing}. This means that we can, for instance, simply pick $G:=-\nabla V:M_c\subseteq X\to T_{M_c}X$.
Since $\partial M_c$ is a compact, orientable, embedded submanifold, we have $H_{n-1}(\partial M_c)\cong \mathbb{Z}^k$ with $k\in \mathbb{N}_{>0}$ being the finite number of disjoint components (\textit{e.g.}, consider $H^0(\partial M_c)$ and Poincar\'e duality, noting that $\mathrm{dim}(\partial M_c)=n-1$). 

Then, for $n\geq 2$ we have
\begin{align*}
    -\nabla V_* H_{n-1}(M_c\setminus A) &\cong -  \nabla V_* H_{n-1}(\partial M_c) \\
    &=  -(\nabla V/\|\nabla V\|_2)_* H_{n-1}(\partial M_c) \cong \mathrm{deg}(h) \cdot \mathbb{Z},
\end{align*}
where $\partial M_c \ni  x\mapsto h(x):=-(\nabla V(x)/\|\nabla V(x)\|_2)\in \mathbb{S}^{n-1}$ and $\mathrm{deg}(h)=\sum_{i=1}^k \mathrm{deg}(h_i)$, for $h_i$ acting on the $i$-th component of $\partial M_c$. Since $\chi(A)\neq 0$, we know from the Poincar\'e-Hopf theorem \cite[p. 35-36]{ref:milnor65} that $\mathrm{deg}(h)\neq 0$. 

For $n=1$, one has $H_0(\mathbb{S}^0)\cong \mathbb{Z}^2$ and hence the degree construction does not readily extend. However, since $-(\nabla V/\|\nabla V\|_2)$ is surjective, one can conclude that $-\nabla V_* H_{n-1}(M_c\setminus A)\cong \mathbb{Z}^2$.

Now, as the spherical adversary condition fails to hold, there is a $z\in \mathbb{S}^{n-1}$ such that $z\neq f(u_x)/\| f(u_x)\|_2$ for all $u_x\in \pi_u^{-1}(N\setminus A)$, for $N$ being some sufficiently small open neighbourhood of $A$.

In what follows we assume that $c>0$ is selected so that $M_c\subseteq N$. This is always possible by compactness of $A$ (and $X$ admitting a metric). 

Now consider the second part of the homology condition. It follows that 
\begin{align*}
    f_* H_{n-1}(\Sigma_{M_c\setminus A}) = (f/\| f\|_2)_*H_{n-1}(\Sigma_{M_c\setminus A}),
\end{align*}
where the normalization is well-defined as $f^{-1}(0)\cap \Sigma_{M_c\setminus A}=\varnothing$. Then, as the normalized map to the sphere is not surjective, consider the factorization
\begin{align*}
    (f/\|f\|_2) : \Sigma_{M_c\setminus A}\overset{\widetilde{f}}{\to}\mathbb{S}^{n-1}\setminus \{*\}\overset{\iota}{\hookrightarrow} \mathbb{S}^{n-1},
\end{align*}
\textit{i.e.}, $f=:\iota \circ \widetilde{f}$, and note that $\mathbb{S}^{n-1}\setminus \{*\}\cong \{*\}$. However, this means that\\ $(f/\| f\|_2)_*H_{n-1}(\Sigma_{M_c\setminus A})$ factors through $H_{n-1}(\{*\})\cong 0$ for $n\geq 2$ and through $H_0(\{*\})\cong \mathbb{Z}$ for $n=1$. Hence, the homology condition fails. 
\end{proof}

We also observe that Theorem~\ref{thm:compare}, \cite[Thm. 4]{ref:kvalheim2023relationships} and the original result by Coron \cite[(B)]{ref:Coron1990}, show that only the top homology is enough to imply the adversary condition, that is, for $\mathrm{dim}(X)=n$, one may consider the ``top homology condition''
\begin{align}
\label{equ:top:homology1}
    G_* H_{n-1} (N\setminus A) \subseteq f_* H_{n-1}(\Sigma_{N\setminus A}) \quad \text{in}\quad  H_{n-1} (T_{N\setminus A}X\setminus 0_{TX}),
\end{align}
for some sufficiently small neighbourhood $N$ of $A$.

\subsection{Generalizations}
\label{sec:generalizaions}
The goal of this note is to clarify the interpretation of Example~\ref{ex:counter} when it comes to comparing the adversary and homology condition. From that perspective, Assumption~\ref{ass:X} is sufficient. However, in general, this assumption is restrictive. First, note that Corollary~\ref{cor:homology} already accommodates general control systems, not relying on Assumption~\ref{ass:X}. Nevertheless, to na\"ively extend the spherical adversary condition {\myref{spherical:adversary:condition}{(iii)}} beyond Assumption~\ref{ass:X}, one may apply again the Poincar\'e-Hopf Theorem and recover the following. 
\\\\
\textbf{\mylabel{sphere:bundle:adversary:condition}{(iiib)} The sphere bundle adversary condition} For any sufficiently small neighbourhood $N\subseteq X$ of $A$ there is local continuous section $\sigma:N\setminus A\to \mathcal{U}$ such that 
\begin{equation}
\label{equ:sphere:bundle:condition}
    \frac{f\circ \sigma}{\| f\circ \sigma \|_g}:N\setminus A\to S(T_{N\setminus A}X)
\end{equation}
is well-defined, but cannot be extended to all of $N$. 
\\\\
Note, this is not precisely a generalization, in Theorem~\ref{thm:spherical:adversary} we can exploit the underlying Euclidean structure and go one step further by exploiting the topological degree\footnote{We remark that there is a recent generalization of Hopf's degree theorem, from spheres to sphere bundles \cite{ref:kvalheim2023hopfgeneralization}.}. What is more, the sense of ``direction'' is lost and the explicit dependency on a continuous section $\sigma$ makes \eqref{equ:sphere:bundle:condition} impractical. This is similar to employing the index condition by Krasnosel'ski\u{\i} directly\footnote{To be more explicit. Consider the setting of Brockett's condition.
Given some vector field $F$ on some $n$-dimensional manifold $X$. Then, any isolated equilibrium point $p$ of $F$ that happens to be locally asymptotically stable will satisfy $\mathrm{ind}_p(F)=(-1)^n$. However, that means that a necessary condition for the existence of a stabilizing feedback $\mu$ is that there is a continuous map $\widetilde{\mu}$, stabilizing or not, such that resulting closed-loop vector field satisfies this index condition.}. Nevertheless, as the spherical adversary condition {\myref{spherical:adversary:condition}{(iii)}}, one may use the sphere bundle adversary condition {\myref{sphere:bundle:adversary:condition}{(iiib)}} to show, for example, that there is no $\mathbb{S}^1$-valued feedback that stabilizes the origin in~\eqref{equ:nonholonomic}. 

Arguably more practical would be a homological condition that does not rely on searching for sections. As alluded to before, following a path that led to adversarial conditions will provide us with a condition on top homology. As only the top degree appears, this shows again the advantage of the homology condition. To be more explicit, as in the proof of Theorem~\ref{thm:spherical:adversary}, given a (local) section $\sigma:N\to \mathcal{U}$, which may be discontinuous on $A$, pick $c>0$ and set $V^{-1}([0,c])=:M_c\subseteq N$. Since $-\nabla  V:\partial M_c\to T_{\partial M_c}X$ is a vector field, we have that $\pi_x\circ -\nabla V = \mathrm{id}_{\partial M_c}$. Hence, 
\begin{align*}
    (\pi \circ -\nabla V)_* H_{n-1}(\partial M_c) = H_{n-1}(\partial M_c) \cong \mathbb{Z}^k,
\end{align*}
for some $k>0$. Thereby, $-\nabla V_* H_{n-1}(\partial M_c) \cong -\nabla V_* H_{n-1}(M_c \setminus A)$ cannot be trivial. Now, we know that $-\nabla V$ is homotopic to $f\circ \sigma$ on $M_c\setminus A$, without introducing new zeros. Hence 
\begin{align*}
     -\nabla V_* H_{n-1}(N \setminus A) = (f\circ \sigma)_* H_{n-1}(N \setminus A)\subseteq f_* H_{n-1}(\Sigma_{N \setminus A})
\end{align*}
cannot be trivial, as a subgroup of $H_{n-1}(T_{N\setminus A}X \setminus 0_{TX})$. Combining the above with Kvalheim's homotopy theorem \cite[Thm. 1]{ref:kvalheim2022obstructions}, we get that
\begin{align*}
\label{equ:top:homology2}
    0\not\cong G_* H_{n-1} (N\setminus A) \subseteq f_* H_{n-1}(\Sigma_{N\setminus A}) \quad \text{in}\quad  H_{n-1} (T_{N\setminus A}X\setminus 0_{TX}).
\end{align*}
Note, deriving this condition on top homology, relies only on the existence of a section $\sigma$ and the non-trivial top homology of $N\setminus A$. Evidently, under mild assumptions, such a condition is implied by the homology condition {\myref{homology:condition}{(ii)}}. Note that the above does \textit{not} rely on $\chi(A)\neq 0$. Hence, when it comes to generality, Corollary~\ref{cor:homology} provides us with the most general necessary condition. However, as adversarial conditions are typically easier to employ, when it comes to applicability, one might hope to improve upon~\eqref{equ:sphere:bundle:condition}. 
 


\begin{remark}[Connection to the work by Mansouri]
\label{rem:Mansouri}
\upshape{
    In \cite{ref:mansouri2007local,ref:mansouri2010topological}, Mansouri provides necessary conditions similar to a top homology condition. For $X=\mathbb{R}^n$, $A$ some $m$-dimensional compact embedded submanifold of $X$ and $f:\Omega\subseteq \mathbb{R}^n\times \mathbb{R}^m\to \mathbb{R}^n$ a $C^1$ map representing the control system, \textit{i.e.}, $\dot{x}=f(x,u)$, Mansouri studies necessary condition for the existence of a $C^1$ feedback such that $A$ is stabilized. Comments on continuous feedback are provided. Specifically, let $\Sigma_{f,\epsilon}:=\{(x,u)\in \Omega:f(x,u)\neq 0,\, d(x,A)<\epsilon\}$, then:
\begin{enumerate}[(a)]
    \item when $m = n-1$ (and $n\geq 2$), Mansouri shows that $\chi(M_A)\cdot \mathbb{Z}\subseteq f_* H_{n-1}(\Sigma_{f,\epsilon})$, for all $\epsilon>0$, where $M_A$ is the manifold enclosed by $A$; on the other hand 
    \item when $m\leq n-2$, $\chi(A)\cdot \mathbb{Z}\subseteq f_* H_{n-1}(\Sigma_{f,\epsilon})$, for all $\epsilon>0$. 
\end{enumerate}
Compared to this work, we observe a difference in the state space, regularity of the attractor and the regularity of the vector fields. Otherwise, one should observe that if one constrains the control systems, sharper necessary conditions can be derived, \textit{e.g.}, we can speak of explicit subgroups instead of only stating that they are non-trivial. Also note the difference between $\Sigma_{f,\epsilon}$ and $\Sigma_{N\setminus A}$. Comparing to {\myref{homology:condition}{(ii)}} and with the above in mind, we believe it is beneficial to consider all degrees. 
}
\exampleEnd
\end{remark}

\section{Discussion and conclusion}
\label{sec:discussion}
\begin{enumerate}[(a)]
    \item (On sufficiency): Clearly, these necessary conditions are not sufficient. As the index condition by Krasnosel'ski\u{\i} and Zabreiko \cite[p. 82]{ref:jongeneel2023topological}, all conditions above fail to capture that $\dot{x}=x$ (sometimes denoted as $\partial_x$), with $x\in \mathbb{R}^{2n}$, cannot be stabilized (\textit{i.e.}, cannot render $0\in \mathbb{R}^{2n}$ asymptotically stable). 
    \item (On $\chi(A)\neq 0$): The standard adversarial condition {\myref{adversary:condition}{(i)}} cannot be extended to the case of $\chi(A)=0$. The same is true for the spherical adversarial condition, as proposed in this note {\myref{spherical:adversary:condition}{(iii)}}. Nevertheless, a top-homological version of {\myref{spherical:adversary:condition}{(iii)}} can be applied, although, trivially. To emphasize the importance of filling this gap, we elaborate on a recurring example by Mansouri, see \cite{ref:mansouri2007local,ref:mansouri2010topological}. Consider the control system
    \begin{equation}
    \label{equ:brockettV2}
        \begin{aligned}
            \dot{x}_1 &= u_1\\
            \dot{x}_2 &= u_2\\
            \dot{x}_3 &= (x_2 u_1 - x_1u_2)e^{{x_3}^2},
        \end{aligned}
    \end{equation}
    with $(x,u)\in \mathbb{R}^3\times \mathbb{R}^2$ (\textit{i.e.}, this control system can be identified with a vector bundle). Note, the system is a minor modification of the one presented in \cite{ref:mansouri2007local}. Now, replacing $u$ by the feedback $x\mapsto (\mu_1(x),\mu_2(x))$ defined through
    \begin{equation}
    \label{equ:mu}
        \begin{aligned}
            \mu_1(x) &:= -x_2x_3 - x_1(x_1^2+x_2^2-1)\\
            \mu_2(x) &:= x_1x_3 - x_2(x_1^2+x_2^2-1)
        \end{aligned}
    \end{equation}
    results in a closed-loop system with $A=\{(x_1,x_2,0)\in \mathbb{R}^3:x_1^2+x_2^2=1\}$ being locally asymptotically stable. See that $(A,X)$ satisfies Assumption~\ref{ass:AX} and $\chi(A)=0$. One may check this using the Lyapunov function $x\mapsto V(x):= (x_1^2+x_2^2-1)^2+x_3^2$. In fact, this Lyapunov function asserts almost global asymptotic stability (\textit{i.e.}, the system is asymptotically stable on $\mathbb{R}^3\setminus \{(0,0,z)\in \mathbb{R}^3:z\in \mathbb{R}\}$).

Control system \eqref{equ:brockettV2} is used in conjunction with feedback \eqref{equ:mu} to show that the adversary condition cannot be readily extended to the case of $\chi(A)=0$ (\textit{i.e.}, the example above admits a stabilizing feedback, but the condition fails).  We argue now that the same is true for the spherical adversary condition.

Indeed, there is no $(x,u)$ such that $f(x,u)=(0,0,\varepsilon)$, for any $\varepsilon>0$, so the spherical adversary condition fails in general, for any compact attractor $A$ such that $\chi(A)\neq 0$. Nonetheless, we can find a stabilizing feedback, where the input is $(\mathbb{S}^1\hookrightarrow \mathbb{R}^2)$-valued and as such, the control system cannot be identified with a vector bundle. In particular, we will project \eqref{equ:mu} onto $\mathbb{S}^1$. Here, a simple computation shows that $\| \mu(x)\|^2_2 = (x_1^2+x_2^2) {V(x)}$, and this expression is only zero on $A\cup \{(0,0,z)\in \mathbb{R}^3:z\in \mathbb{R}\}$. Using \eqref{equ:brockettV2} and the Lyapunov function $V$ from above, one observes that normalizing $\mu$, preserves asymptotic stability (now understood in the sense of Filippov). Moreover, under the normalized feedback
    \begin{align*}
        \langle \nabla V(x),\dot{x}\rangle &= -\frac{4}{\|\mu(x)\|_2}(x_1^2+x_2^2)(x_1^2+x_2^2-1) - \frac{2}{\|\mu(x)\|_2}(x_1^2+x_2^2)x_3^2 e^{x_3^2}\\
         &\leq  -\frac{4}{\|\mu(x)\|_2}(x_1^2+x_2^2)(x_1^2+x_2^2-1) - \frac{2}{\|\mu(x)\|_2}(x_1^2+x_2^2)x_3^2,
    \end{align*}
    as $\inf_{z\in \mathbb{R}} \{e^{z^2}\}=1$.
    Then, recall that for any $c>0$, $M_c=V^{-1}([0,c])$ is compact and forward invariant. Now, pick any $c\in (0,1)$ and see that on $M_c$ we have $(x^2+y^2-1)^2< 1$ and thus $(x^2+y^2-1)^2 < (x^2+y^2-1)$. What is more, for such a $c>0$, see that $\min_{x\in M_c} \{x_1^2 + x_2^2\} = 1-\sqrt{c}$. Hence, we can write
    \begin{align*}
        \langle \nabla V(x),\dot{x}\rangle\big|_{M_c} &\leq  -\frac{4}{\|\mu(x)\|_2}(x_1^2+x_2^2)(x_1^2+x_2^2-1)^2 - \frac{2}{\|\mu(x)\|_2}(x_1^2+x_2^2)x_3^2\\
        &\leq -\frac{1}{\|\mu(x)\|_2}(x_1^2+x_2^2)V(x)\\
        &= - \sqrt{(x_1^2+x_2^2)}\sqrt{V(x)}\\
        &\leq  - \sqrt{1-\sqrt{c}}\sqrt{V(x)}=:-\alpha \sqrt{V}(x),
    \end{align*}
    for $\alpha\in (0,1)$. This means, however, that $A$ is locally finite-time asymptotically stable (\textit{e.g.}, on $M_c$) under \eqref{equ:brockettV2} and the feedback $x\mapsto \mu(x)/\| \mu(x)\|_2$ (with $\mu(x):=0$ on $A$), \textit{e.g.}, see \cite[Thm. 4.2]{ref:bhat2000finitetime}. This, despite not satisfying the spherical adversary condition {\myref{spherical:adversary:condition}{(iii)}}. 
    
    A top-homological generalization, and in particular the homology condition {\myref{homology:condition}{(ii)}}, do not rely on $\chi(A)=0$ and must be true. However, as in Remark~\ref{rem:Mansouri}, several conditions become trivial. We leave the study of non-trivial necessary conditions for circular attractor for future work, see also \cite[p. 538]{ref:kvalheim2022obstructions}.

    \item (On Control Barrier Functions (CBFs)): Consider the control system $\dot{x}=x+u$, with $x\in \mathbb{R}^2$ and $u\in \mathbb{S}^1\hookrightarrow \mathbb{R}^2$. In \cite{ref:de2026topological}, we find Brockett-like necessary condition for the existence of a continuous controller such that some set $\mathcal{S}$ is safe (forward invariant). They consider $\mathcal{S}:=\{x\in \mathbb{R}^2: \|x\|_2^2\leq 1\}$ and show that there is no continuous feedback \cite[Ex. 4]{ref:de2026topological}, employing effectively the same arguments as in Example~\ref{ex:counter} (\cite[Ex. 4]{ref:kvalheim2023relationships}). However, one may select the controller $\kappa:\mathbb{R}^2\to \mathbb{S}^1\hookrightarrow \mathbb{R}^2$, as defined through $x\mapsto \kappa(x):=-x/\|x\|_2$ for $x\in \mathbb{R}^2\setminus \{0\}$ while $\kappa(x):=0$ otherwise. Indeed, as the closed-loop system becomes
    \begin{equation*}
        \dot{x} = \begin{cases}
            \left(\dfrac{\|x\|_2-1}{\|x\|_2}\right)x,\quad &\text{for}\,\, x\in \mathbb{R}^2\setminus \{0\}\\
            0\quad &\text{otherwise},
        \end{cases}
    \end{equation*}
    we render $\mathcal{S}$ safe (and $0\in \mathcal{S}$ is asymptotically stable in the sense of Filippov). 

    \item (On discontinuous feedback): We emphasize that the type of discontinuities we consider are mild (see \S\ref{sec:discont:control:system}) and are not intended to overcome topological obstructions. Indeed, typically, if there is no continuous feedback, then, there is no \textit{robust} discontinuous feedback either \cite{ref:ryan1994brockett,ref:ceragioli2002some,ref:mayhew2011topological}. For us, discontinuities emerge due to the possibility of working with fiber bundles with non-contractible fibers \textit{cf}. Example~\ref{ex:counter}-\ref{ex:finite:time:sphere} and the robustness enters due to application of the Filippov operator. Handling discontinuities outside of $A$ is left for future work, see for instance \cite[Ex. 3]{ref:ceragioli2002some} and \cite{ref:montgomery2024hybrid}.     

\end{enumerate}
    In this work, the study of certain control systems is enabled by the existence of a smooth Lyapunov function. In a line of work pioneered by Lyapunov, Morse, Conley and many others, this note may be interpreted as another example where classification of dynamical (control) systems through the existence of these smooth functions is fruitful. One should compare this with classification of systems based on the regularity of the vector field itself. 

    Considering more general control systems, \textit{e.g.}, allowing more general discontinuities and fiber bundles, is left for future work.

\section*{Acknowledgments}
The first author would like to thank Matthew D. Kvalheim for helpful discussions.

\bibliographystyle{siamplain}
\bibliography{references}
\end{document}

%% file: ex_shared.tex
\usepackage{lipsum}
\usepackage{amsfonts}
\usepackage{graphicx}
\usepackage{epstopdf}
\usepackage{algorithmic}
\ifpdf
  \DeclareGraphicsExtensions{.eps,.pdf,.png,.jpg}
\else
  \DeclareGraphicsExtensions{.eps}
\fi

\newsiamremark{remark}{Remark}
\newsiamremark{hypothesis}{Hypothesis}
\crefname{hypothesis}{Hypothesis}{Hypotheses}
\newsiamthm{claim}{Claim}
\newsiamremark{fact}{Fact}
\crefname{fact}{Fact}{Facts}

\headers{On necessary conditions for feedback stabilizability}{W. Jongeneel, M. Scolamiero and K. H. Johansson}

\title{On necessary conditions
for feedback stabilizability when the control space is\\ a fiber bundle\thanks{\funding{This work received funding from Digital Futures through the project ``\textit{Programmability of Cells}".}}}
\author{Wouter Jongeneel \footnotemark[2]\,\,\footnotemark[3] \and Martina Scolamiero \thanks{Department of Mathematics, KTH Royal Institute of Technology, Stockholm, Sweden 
  (\email{wouterjo@kth.se}, \email{scola@kth.se}).} \and Karl H. Johansson \thanks{Department of Decision and Control Systems, KTH Royal Institute of Technology, Stockholm, Sweden (\email{wouterjo@kth.se}, \email{kallej@kth.se}).}}

\usepackage{amsopn}


%% file: references.bib
@book{ref:munkrestopology,
  title={Topology},
  author={Munkres, James},
  year={2014},
publisher={Pearson},
address = {Harlow}
}

@article{ref:kvalheim2023hopfgeneralization,
  title={A generalization of the {H}opf degree theorem},
  author={Kvalheim, Matthew D.},
  journal={Proc. Am. Math. Soc.},
  volume={151},
  number={01},
  pages={453--454},
  year={2023}
}

@book{ref:jongeneel2023topological,
  title={Topological Obstructions to Stability and Stabilization: History, Recent Advances and Open Problems},
  author={Jongeneel, Wouter and Moulay, Emmanuel},
  year={2023},
  publisher={Springer Nature},
address={Cham}
}

@article{ref:Coron1990,
title = {A necessary condition for feedback stabilization},
journal = {Syst. Control Lett.},
volume = {14},
number = {3},
pages = {227--232},
year = {1990},
author = {Jean-Michel Coron}
}

@article{ref:wilson1967structure,
  title={The structure of the level surfaces of a {L}yapunov function},
  author={Wilson Jr, F Wesley},
  journal={J. Differ. Equ.},
  volume={3},
  number={3},
  pages={323--329},
  year={1967}
}

@INPROCEEDINGS{ref:brockett1983asymptotic,
    author = {R. W. Brockett},
    title = {Asymptotic stability and feedback stabilization},
    booktitle = {Differential Geometric Control Theory},
    year = {1983},
    pages = {181--191},
    publisher = {Birkh\"auser},
    address={Boston}
}

@article{ref:clarke1998asymptotic,
  title={Asymptotic stability and smooth {L}yapunov functions},
  author={Clarke, Francis H and Ledyaev, Yu S and Stern, Ronald J},
  journal={J. Differ. Equ.},
  volume={149},
  number={1},
  pages={69--114},
  year={1998}
}

@article{ref:moulay2011conley,
  title={Conley index condition for asymptotic stability},
  author={Moulay, Emmanuel and Hui, Qing},
  journal={Nonlinear Anal.-Theor.},
  volume={74},
  number={13},
  pages={4503--4510},
  year={2011},
  publisher={Elsevier}
}

@INPROCEEDINGS{ref:jongeneelECC24,
  author={Jongeneel, Wouter},
  booktitle={2024 European Control Conference (ECC)}, 
  title={On topological properties of compact attractors on {H}ausdorff spaces}, 
  year={2024},
  volume={},
  number={},
  pages={186-191}}

@article{ref:mayhew2011topological,
  title={On the topological structure of attraction basins for differential inclusions},
  author={Mayhew, Christopher G and Teel, Andrew R},
  journal={Syst. Control Lett.},
  volume={60},
  number={12},
  pages={1045--1050},
  year={2011}
}

@book{ref:milnor65,
  title={Topology from the Differentiable Viewpoint},
  author={Milnor, John},
  year={1965},
  publisher={Princeton University Press},
  address={Princeton}
}

@book{ref:bredon1993topology,
  title={Topology and geometry},
  author={Bredon, Glen E},
  year={1993},
  publisher={Springer},
  address={New York}
}

@book{ref:krasnosel1968operator,
  title={The operator of translation along the trajectories of differential equations},
  author={Krasnosel'ski\u{\i}, Mark A},
  year={1968},
  publisher={American Mathematical Society},
  address={Providence}
}

@incollection{ref:montgomery2024hybrid,
  author    = {Montgomery, R. and Sanfelice, R.},
  title     = {Hybrid Control, {M}orse Theory and {I}van {K}upka},
  booktitle = {Ivan Kupka Legacy: A Tour Through Controlled Dynamics},
  editor    = {Bernard Bonnard and Monique Chyba and David Holcman and Emmanuel Tr\'elat},
  publisher = {American Institute of Mathematical Sciences},
  address = {Springfield},
  year      = {2024},
  pages     = {135--153}
}

@article{ref:bhat2000finitetime,
  title={Finite-time stability of continuous autonomous systems},
  author={Bhat, Sanjay P and Bernstein, Dennis S},
  journal={SIAM J. Contr. Optim.},
  volume={38},
  number={3},
  pages={751--766},
  year={2000},
  publisher={SIAM}
}

@article{ref:de2026topological,
  title={Topological Obstructions to the Existence of Control Barrier Functions},
  author={de Sa, Massimiliano and Ames, Aaron D},
        journal = {arXiv e-prints},
         year = 2026,
       eprint = {2603.18422},
}

@incollection{ref:kvalheim2023relationships,
  author    = {Kvalheim, Matthew D},
  title     = {Relationships between necessary conditions for feedback stabilizability},
  booktitle = {Geometry, Topology and Control System Design: Proceedings of a Banff International Research Station Workshop},
  editor    = {Mohamed Ali Belabbas},
  publisher = {American Institute of Mathematical Sciences},
  address = {Springfield},
  year      = {2025},
  pages     = {167--179}
}

@incollection{ref:whydowetwist2025,
  author    = {Belabbas, Mohamed-Ali and Ko, Jehyung},
  title     = {Dynamic Global Feedback Stabilization: why do the twist?},
  booktitle = {Geometry, Topology and Control System Design: Proceedings of a Banff International Research Station Workshop},
  editor    = {Mohamed Ali Belabbas},
  publisher = {American Institute of Mathematical Sciences},
  address = {Springfield},
  year      = {2025},
  pages     = {47--60}
}

@incollection{ref:Mansouri2025,
  author    = {Abdol-Reza Mansouri
},
  title     = {On {C}oron’s Necessary Condition for Local Asymptotic
Stabilization},
  booktitle = {Geometry, Topology and Control System Design: Proceedings of a Banff International Research Station Workshop},
  editor    = {Mohamed Ali Belabbas},
  publisher = {American Institute of Mathematical Sciences},
  address = {Springfield},
  year      = {2025},
  pages     = {337--340}
}

@article{ref:cai2008smooth2,
  title={Smooth {L}yapunov functions for hybrid systems part {II}: (pre) asymptotically stable compact sets},
  author={Cai, Chaohong and Teel, Andrew R and Goebel, Rafal},
  journal={IEEE T. Automat. Control},
  volume={53},
  number={3},
  pages={734--748},
  year={2008},
  publisher={IEEE}
}

@book{ref:coron2007control,
  title={Control and Nonlinearity},
  author={Coron, Jean-Michel},
  year={2007},
  publisher={American Mathematical Society},
  address={Providence}
}

@article{ref:kvalheim2021necessary,
  title={Necessary conditions for feedback stabilization and safety},
  author={Kvalheim, Matthew D and Koditschek, Daniel E},
        journal = {J. Geom. Mech.},
         volume = {14},
number = {4},
pages = {659-693},
year = {2022}
}

@article{ref:jongeneel2025topological,
  title={On topological properties of closed attractors},
  author={Jongeneel, Wouter},
       eprint={2511.10429},
      archivePrefix={arXiv},
  year={2025}
}

@article{ref:kvalheim2026stabilizability,
  title={Stabilizability of first-order dynamics in second-order systems},
  author={Kvalheim, Matthew D},
eprint={2604.10757},
        archivePrefix={arXiv},
  year={2026}
}

@article{ref:byrnes2008brockett,
  title={On {B}rockett's Necessary Condition for Stabilizability and the Topology of {L}iapunov Functions on $\mathbb{R}^n$},
  author={Byrnes, Christopher I},
  journal={Commun. Inf. Syst.},
  volume={8},
  number={4},
  pages={333--352},
  year={2008}
}

@book{Hatcher,
  title={Algebraic Topology},
  author={Hatcher, Allen},
  year={2002},
  publisher={Cambridge University Press},
  address={Cambridge}
}

@article{ref:gobbino2001topological,
  title={Topological properties of attractors for dynamical systems},
  author={Gobbino, Massimo},
  journal={Topology},
  volume={40},
  number={2},
  pages={279--298},
  year={2001},
  publisher={Elsevier}
}

@article{ref:kapitanski2000shape,
  title={Shape and {M}orse theory of attractors},
  author={Kapitanski, Lev and Rodnianski, Igor},
  journal={Comm. Pure Appl. Math.},
  volume={53},
  number={2},
  pages={218--242},
  year={2000},
  publisher={Wiley Online Library}
}

@book{ref:krasnosel1984geometrical,
  title={Geometrical methods of nonlinear analysis},
  author={Krasnosel'ski\u{\i}, Alexander and Zabreiko, Petr Petrovich},
  year={1984},
  publisher={Springer},
  address={Berlin}
}

@article{ref:bobylev1974deformation,
  title={Deformation of a system into an asymptotically stable system},
  author={Bobylev, N.A. and Krasnosel'ski\u{\i}, M.A.},
  journal={Autom. Remote Control},
  volume={35},
  number={7},
  pages={1041--1044},
  year={1974}
}

@book{ref:guillemin2010differential,
  title={Differential Topology},
  author={Guillemin, Victor and Pollack, Alan},
  year={2010},
  publisher={American Mathematical Society},
  address={Providence}
}

@book{ref:wonham1974linear,
  title={Linear multivariable control},
  author={Wonham, W Murray},
  year={1979},
  publisher={Springer},
  address={New York}
}

@book{ref:bhatia1970stability,
  title={Stability theory of dynamical systems},
  author={Bhatia, Nam Parshad and Szeg{\"o}, Giorgio P},
  year={1970},
  publisher={Springer},
  address={Berlin}
}

@article{ref:mansouri2010topological,
  title={Topological obstructions to submanifold stabilization},
  author={Mansouri, Abdol-Reza},
  journal={IEEE T. Automat. Contr.},
  volume={55},
  number={7},
  pages={1701--1703},
  year={2010}
}

@article{ref:mansouri2007local,
  title={Local asymptotic feedback stabilization to a submanifold: Topological conditions},
  author={Mansouri, Abdol-Reza},
  journal={Syst. Control Lett.},
  volume={56},
  number={7-8},
  pages={525--528},
  year={2007}
}

@article{ref:fathi2019smoothing,
  title={Smoothing {L}yapunov functions},
  author={Fathi, Albert and Pageault, Pierre},
  journal={T. Am. Math. Soc.},
  volume={371},
  number={3},
  pages={1677--1700},
  year={2019}
}

@inproceedings{ref:sontag1980remarks,
  title={Remarks on continuous feedback},
  author={Sontag, Eduardo and Sussmann, Hector},
  booktitle={Proc. IEEE Conf. Decis. Control},
  pages={916--921},
  year={1980}
}

@article{ref:kvalheim2022obstructions,
author = {Kvalheim, Matthew D.},
title = {Obstructions to Asymptotic Stabilization},
journal = {SIAM J. Contr. Optim.},
volume = {61},
number = {2},
pages = {536-542},
year = {2023}
}

@article{ref:zabczyk1989,
  title={Some comments on stabilizability},
  author={Zabczyk, Jerzy},
  journal={Appl. Math. Opt.},
  volume={19},
  number={1},
  pages={1--9},
  year={1989}
}

@book{ref:nonlinear1990NvdS,
   author = "H. Nijmeijer and A.J. van der Schaft",
   publisher = "Springer",
   title = "Nonlinear Dynamical Control Systems",
   year = "1990",
   address={New York}
}

@book{ref:HatcherVBK2003,
  title={Vector bundles and {K}-theory},
  author={Hatcher, Allen},
  url={http://www. math. cornell. edu/~hatcher},
  year={2003}
}

@article{ref:ceragioli2002some,
  title={Some remarks on stabilization by means of discontinuous feedbacks},
  author={Ceragioli, Francesca},
  journal={Syst. Control Lett.},
  volume={45},
  number={4},
  pages={271--281},
  year={2002}
}

@article{ref:ryan1994brockett,
  title={On {B}rockett’s condition for smooth stabilizability and its necessity in a context of nonsmooth feedback},
  author={Ryan, Eugene P},
  journal={SIAM J. Contr. Optim.},
  volume={32},
  number={6},
  pages={1597--1604},
  year={1994}
}

@article{ref:orsi2003necessary,
  title={Necessary conditions for stability and attractivity of continuous systems},
  author={Orsi, Robert and Praly, L and Mareels, Iven},
  journal={Int. J. Control},
  volume={76},
  number={11},
  pages={1070--1077},
  year={2003}
}

@book{Lee2,
  title={Introduction to Smooth Manifolds},
  author={Lee, John M.},
  year={2012},
  publisher={Springer},
  address={New York}
}

@book{ref:filippov1988differential,
  title={Differential Equations with Discontinuous Righthand Sides},
  author={Filippov, A F},
  year={1988},
  publisher={Springer Science \& Business Media},
  address={Dordrecht}
}
